# Projection pursuit for discrete data

**Persi Diaconis[1] and Julia Salzman[*2]**

*Stanford University*

**Abstract:** This paper develops projection pursuit for discrete data using the discrete Radon transform. Discrete projection pursuit is presented as an exploratory method for finding informative low dimensional views of data such as binary vectors, rankings, phylogenetic trees or graphs. We show that for most data sets, most projections are close to uniform. Thus, informative summaries are ones deviating from uniformity. Syllabic data from several of Plato's great works is used to illustrate the methods. Along with some basic distribution theory, an automated procedure for computing informative projections is introduced.

## 1. Introduction

Projection pursuit is an exploratory graphical tool for picturing high dimensional data through low dimensional projections. Introduced by Kruskal [35], [36], and developed by Friedman and Tukey [28], the idea is to have the computer select a small family of projections by numerically optimizing an index of "interest". The original projection indices were ad hoc. In joint work with David Freedman [15], it was shown that for most data sets, most projections are about the same: approximately Gaussian. Therefore, the interesting projections, the ones which were special for this data set, are projections that are far from Gaussian.

Peter Huber [32] found his own version of this: projections are uninformative if they are unstructured or "random". Thus projections with high entropy are uninformative. For a fixed scale, distributions having high entropy are approximately Gaussian. Huber also showed that the Friedman-Tukey index is a measure of non-Gaussianess.

The purpose of the present paper is to give a parallel development for data in discrete spaces: collections of binary vectors, rankings, phylogenetic trees or sets of graphs. We develop a notion of projection as a partition of the discrete data into blocks. We show that most for most data sets, most projections are close to uniformly partitioned. This suggests that the informative summaries are the ones with splits that are far from uniform.

The outline of the paper is as follows. Definitions and first examples are given in Section 2. The ideas lean on classical developments in block designs and give new applications for that theory. A discrete version of the Radon transform along with an inversion theory is presented, determining when a collection of projections

---

*Supported in part by NSF Grant DMS-02-41246.

[1]Stanford University, Department of Mathematics and Department of Statistics, Sequoia Hall, 390 Serra Mall, Stanford CA, 94305, USA

[2]Stanford University, Department of Statistics, Sequoia Hall, 350 Serra Mall, Stanford CA 94305, USA, e-mail: julia.salzman@gmail.com

*AMS 2000 subject classifications:* 44A12, 62K10, 90C08.

*Keywords and phrases:* binary vector, discrete data, discrete Radon transform, least uniform partition, phylogenetic tree, projection pursuit, ranking, syllable patterns.

loses information. Section 3 gives a data analytic example in some detail. The data arise from the problem of putting some of Plato's works in chronological order. Here, discrete projection pursuit leads to the discovery of a striking, easily interpretable structure that does not appear in other analyses of this data (eg. Ahn et al. [1], Cox and Brandwood [11], Charnomordic and Holmes [8], Wishart and Leach [49]). Section 4 proves that for most data sets, most partitions lead to approximately uniform projections. This leads directly to a usable criteria: a projection is interesting if it is far from uniform. The distance to uniformity can be measured by any distance between probabilities, and we consider the well-known total variation, Hellinger and Vasserstein metrics.

The final section gives results for the least uniform projection. Theorem 5.1 shows that if the class of projections is not too rich, for example, the affine hyperplane in $\mathbb{Z}_k$, then for most data sets even the least uniform partition is close to uniform. If the class of projections contains many sets, then least uniform projections are "structured". The final theorem attacks the problem of a data analyst finding "structure" in "noise". Computational details for computing the metrics and automating the analysis are in an Appendix.

There has been extensive development of projection pursuit for density estimation (Friedman et al. [26]), regression (Friedman and Stuetzle [27], Hall [30]), applications to time series (Donoho [22]), discriminant analysis (Posse [42], Polzehl [41]) and standard multivariate problems such as covariance estimation (Hwang et al. [33]). This has led to a healthy development captured in the modern implementations (Xgobi, Ggobi). Online documentation for this software is an instructive catalog. We have not attempted to develop our ideas in these directions, but the beginning steps of ridge functions will be found below.

## 2. Projections and Radon transforms

This section introduces our notation and set up for working with discrete data. It defines projection bases, the discrete Radon transform and gives examples with binary data and permutation data. Analysis will be performed on binary $n$-tuple data from several works of Plato. Let $\mathcal{X}$ be a finite set. Let $\mathcal{Y}$ be a class of subsets of $\mathcal{X}$. Let $f : \mathcal{X} \to \mathbb{R}$ be a function. The Radon transform of $f$ at $y \in \mathcal{Y}$ is defined by

$$\bar{f}(y) = \sum_{x \in y} f(x). \tag{1}$$

The class $\mathcal{Y}$ is called a projection base if:

$$|y| \text{ is constant for } y \in \mathcal{Y} \quad (|y| \text{ denotes the cardinality of } \mathcal{Y}). \tag{2}$$

(3)     There is a partition $p_1, \ldots, p_j$ of $\mathcal{Y}$ such that each $p_i$ is a partition of $\mathcal{X}$.

For a partition $p$, the numbers $\bar{f}(y)_{y \in p}$ will be called the projection of $f$ in direction $p$. The sets in $\mathcal{Y}$ may be thought of as "lines" in a geometry. If lines in the same partition are called parallel, then (3) corresponds to the Euclidean axiom: for every point $x \in \mathcal{X}$ and every line $y \in \mathcal{Y}$, there is a unique line $y^*$ parallel to $y$ such that $x \in y^*$. In the statistics literature, designs with property (3) are called "resolvable" (See Hedayat et al. [31] or Constantine [10] for examples). Assumption (2) guarantees that projections are based on averages over comparable sets.

Consider the following examples:

TABLE 1
*Percentage distribution of sentence endings*

| Type of ending | Rep. | Laws | Phil. | Pol. | Soph. | Tim. |
|---|---|---|---|---|---|---|
| ∪ ∪ ∪ ∪ ∪ | 1.1 | 2.4 | 2.5 | 1.7 | 2.8 | 2.4 |
| - ∪ ∪ ∪ ∪ | 1.6 | 3.8 | 2.8 | 2.5 | 3.6 | 3.9 |
| ∪ - ∪ ∪ ∪ | 1.7 | 1.9 | 2.1 | 3.1 | 3.4 | 6.0 |
| ∪ ∪ - ∪ ∪ | 1.9 | 2.6 | 2.6 | 2.6 | 2.6 | 1.8 |
| ∪ ∪ ∪ - ∪ | 2.1 | 3.0 | 4.0 | 3.3 | 2.4 | 3.4 |
| ∪ ∪ ∪ ∪ - | 2.0 | 3.8 | 4.8 | 2.9 | 2.5 | 3.5 |
| - - ∪ ∪ ∪ | 2.1 | 2.7 | 4.3 | 3.3 | 3.3 | 3.4 |
| - ∪ - ∪ ∪ | 2.1 | 1.8 | 1.5 | 2.3 | 4.0 | 3.4 |
| - ∪ ∪ - ∪ | 2.8 | 0.6 | 0.7 | 0.4 | 2.1 | 1.7 |
| - ∪ ∪ ∪ - | 4.6 | 8.8 | 6.5 | 4.0 | 2.3 | 3.3 |
| ∪ - - ∪ ∪ | 3.3 | 3.4 | 6.7 | 5.3 | 3.3 | 3.4 |
| ∪ - ∪ - ∪ | 2.6 | 1.0 | 0.6 | 0.9 | 1.6 | 3.2 |
| ∪ - ∪ ∪ - | 4.6 | 1.1 | 0.7 | 1.0 | 3.0 | 2.7 |
| ∪ ∪ - - ∪ | 2.6 | 1.5 | 3.1 | 3.1 | 3.0 | 3.0 |
| ∪ ∪ - ∪ - | 4.4 | 3.0 | 1.9 | 3.0 | 3.0 | 2.2 |
| ∪ ∪ ∪ - - | 2.5 | 5.7 | 5.4 | 4.4 | 5.1 | 3.9 |
| - - - ∪ ∪ | 2.9 | 4.2 | 5.5 | 6.9 | 5.2 | 3.0 |
| - - ∪ - ∪ | 3.0 | 1.4 | 0.7 | 2.7 | 2.6 | 3.3 |
| - - ∪ ∪ - | 3.4 | 1.0 | 0.4 | 0.7 | 2.3 | 3.3 |
| - ∪ - - ∪ | 2.0 | 2.3 | 1.2 | 3.4 | 3.7 | 3.3 |
| - ∪ - ∪ - | 6.4 | 2.4 | 2.8 | 1.8 | 2.1 | 3.0 |
| - ∪ ∪ - - | 4.2 | 0.6 | 0.7 | 0.8 | 3.0 | 2.8 |
| ∪ ∪ - - - | 2.8 | 2.9 | 2.6 | 4.6 | 3.4 | 3.0 |
| ∪ - ∪ - - | 4.2 | 1.2 | 1.3 | 1.0 | 1.3 | 3.3 |
| ∪ - - ∪ - | 4.8 | 8.2 | 5.3 | 4.5 | 4.6 | 3.0 |
| ∪ - - - ∪ | 2.4 | 1.9 | 5.3 | 2.5 | 2.5 | 2.2 |
| ∪ - - - - | 3.5 | 4.1 | 3.3 | 3.8 | 2.9 | 2.4 |
| - ∪ - - - | 4.0 | 3.7 | 3.3 | 4.9 | 3.5 | 3.0 |
| - - ∪ - - | 4.1 | 2.1 | 2.3 | 2.1 | 4.1 | 6.4 |
| - - - ∪ - | 4.1 | 8.8 | 9.0 | 6.8 | 4.7 | 3.8 |
| - - - - ∪ | 2.0 | 3.0 | 2.9 | 2.9 | 2.6 | 2.2 |
| - - - - - | 4.2 | 5.2 | 4.0 | 4.9 | 3.4 | 1.8 |
| no. sentences | 3778 | 3783 | 958 | 770 | 919 | 762 |

**Example 2.1.** $\mathcal{X} = \mathbb{Z}_2^k$ the set of binary $k$-tuples. Here is a concrete example of a data set with this structure; L. Brandwood classified each sentence of Plato's *Republic* according to its last five syllables. These can run from all short (∪) through all long (-). Identifying ∪ with 1 and - with 0, each sentence is associated with a binary 5-tuple. As $x$ ranges over $\mathbb{Z}_2^5$, let $f(x)$ denote the proportion of sentences with ending $x$. The values of $f(x)$ are given in the first column of Table 1.

A second example of data with this structure is the result of grading correct/incorrect in a test with $k$ questions. There are several useful choices of $\mathcal{Y}$ given next:

### 2.1. Projections for data in $\mathbb{Z}_2^k$

#### 2.1.1. Marginal projections in $\mathbb{Z}_2^k$

For $i = 1, 2, \ldots, k$, let $y_i^0 = \{x \in \mathbb{Z}_2^k : x_i = 0\}$, let $y_i^1 = \{x \in \mathbb{Z}_2^k : x_i = 1\}$. The sets $\mathcal{Y} = \{y_i^j\}, 1 \leq i \leq k, j \in \{0, 1\}$ form a projection base. In the Plato example, the projections have a simple interpretation as the proportion of sentences with a specific syllable in the $i^{\text{th}}$ place. Displaying projections offers no problem here; a single number suffices.

A second natural choice of $\mathcal{Y}$ gives second order margins. This is based on sets $y_{ij}^{ab} = \{x \in \mathbb{Z}_2^k : x_i = a, x_j = b\}, 1 \leq i < j \leq k, \quad a, b \in \{0, 1\}$. In this case, a projection consists of 4 numbers. In the Plato example, the projection along coordinates $i, j$ gives the proportion of sentences with each of the 4 possibles patterns $\cup \cup, \cup$ -, - $\cup$, - - in positions $i, j$. Table 3 in Section 3 is an example of one method to display such projections. Section 2 contains an analysis of the data in Table 1 based on these projections. The analysis gives a clear interpretation to a classical way of dating the books of Plato. The analysis is independent of the other examples in this section and can be read at this time.

Here are some examples to show how the structure of $f$ is reflected in $\bar{f}$. If $f(x) = \delta_{x, x_0}, \bar{f}(y) = 1$ if $x_0 \in y$ and zero otherwise. If $f(x) = \frac{1}{2^k}$ for all $x$, then $\bar{f}(y) = \frac{|y|}{2^k}$ and hence is constant for all $y$. As a final example, consider a fixed, non-zero vector $y^* \in \mathbb{Z}_2^k$. Let $S$ be the hyperplanes determined by $y^* : S = \{x \in \mathbb{Z}_2^k : x \cdot y = 0 \bmod 2\}$. Let

$$f(x) = \begin{cases} \frac{1}{2^k}, & \text{if } x \in S, \\ 0, & \text{otherwise.} \end{cases}$$

An easy computation shows

$$\bar{f}(y_z^0) = \begin{cases} 1, & \text{if } z = y^*, \\ \frac{1}{2}, & \text{otherwise,} \end{cases}$$

$$\bar{f}(y_z^1) = \begin{cases} 0, & \text{if } z = y^*, \\ \frac{1}{2}, & \text{otherwise.} \end{cases}$$

The hyperplane transform is essentially the same as the ordinary Fourier transform on the group $\mathbb{Z}_2^k$. This is defined by

$$\hat{f}(z) = \sum_x (-1)^{x \cdot z} f(x).$$

If $f$ is a probability on $\mathbb{Z}_2^k$, $\hat{f}(z) = 2\bar{f}(y_z^0) - 1$. The transform $\hat{f}$ has been widely used for data analysis of this type of data. See Solomon [44] or Diaconis [18, 19], Chapter 11. The discrete Radon transform with projections onto affine hyperplanes is also used by Ahn et al. [1].

### 2.1.2. *Affine hyperplanes in $\mathbb{Z}_2^k$*

This is one natural way of "filling out" the marginal projections presented above. For $z \in \mathbb{Z}_2^k$ and $a \in \{0, 1\}$, let $y_z^a = \{x \in \mathbb{Z}_2^k : x \cdot z = a \bmod 2\}$. The collection $\mathcal{Y} = \{y_z^a\}_{z \in \mathbb{Z}_2^k, a \in \{0, 1\}}$ forms a projection base. Observe that when $z$ has a 1 in position $i$ and zeros elsewhere, $y_z^a$ equals the $y_i^a$ of the previous example. The sets in $\mathcal{Y}$ are the affine hyperplanes in $\mathbb{Z}_2^k$. Similarly, the affine planes of any dimension form a projection base. An analysis of the Plato data using all affine hyperplanes is in Appendix A.3 below.

## 2.2. **Projections for data in $\mathcal{X} = S_n$, the sets of permutations of $n$ letters.**

Permutation data arises in taste testing, ranking and elections; for example, in presidential elections of the American Psychological Association, members are asked to rank order 5 candidates. Here, for a permutation $\pi$, $f(\pi)$ is taken as the proportion of voters choosing the order $\pi$. For background and many examples, see Critchlow [12], Fligner and Verducci [25] or Marden [39].

### 2.2.1. *Partitions based on marginal projections of permutations in $S_n$.*

Let $y_{ij} = \{\pi \in S_n : \pi(i) = j, \, 1 \leq i, j \leq n\}$. These sets form a projection base. For fixed $i$, the sets $y_{i1}, y_{i2}, \ldots, y_{in}$ form a partition $p(i)$. The projection in direction $p(i)$ has a natural interpretation in the example: how did people rank candidate $i$? The projection can be displayed by making a histogram.

A second useful choice of $\mathcal{Y}$ is based on considering two positions: $y_{ij}^{kl} = \{\pi \in S_n : \pi(i) = k, \pi(j) = l\}$ with $i \neq j, k \neq l$. This leads to projections giving the joint rankings of a fixed pair of candidates in the example. Such projections can be displayed by making a 2-dimensional picture and gray scaling the $(i, j)$ square to correspond to the proportion of voters ranking the pair of candidates in order $(i, j)$. Similarly third and higher order projections can be defined.

### 2.2.2. *Partitions based on subgroups of $S_n$.*

When $\mathcal{X}$ is a group such as $S_n$, the following constructions for $\mathcal{Y}$ are available. Let $N$ be a subgroup of $\mathcal{X}$. The orbits of $N$ acting on $\mathcal{X}$ are the cosets $\{Ny\}_{y \in \mathcal{X}}$, and the distinct orbits partition $\mathcal{X}$. Varying $N$ by conjugation, $\{yNy^{-1}\}_{y \in \mathcal{X}}$, gives a projection base for $\mathcal{X}$.

When $N$ is taken as $S_{n-1} = \{\pi \in S_n : \pi(1) = 1\}$ the projections are the marginal projections defined above. Taking $N$ as $S_{n-2} = \{\pi \in S_n : \pi(1) = 1, \pi(2) = 2\}$ gives the second order margins. An important class of subgroups are the so-called Young subgroups: let $\lambda_1 \leq \lambda_2 \leq \cdots \leq \lambda_n$ be a partition of $n$ so $\sum_i \lambda_i = n$. Let $S_{\lambda_1} \times S_{\lambda_2} \times \cdots \times S_{\lambda_n}$ be the permutations that permute the first $\lambda_1$ elements among themselves and the next $\lambda_2$ elements among themselves, etc. These include the previous examples and provide enough transforms for an inversion theory, as will be shown below. Display of such projections is not a well studied problem. In the case of a projection corresponding to a Young subgroup, one suggestion is a 1-dimensional histogram using one of the orderings suggested in Chapter 3 of James [34].

If $\mathcal{X} = G/H$ where $G$ is a group and $H$ is a subgroup and $G \subset N \subset H$, with $N$ a subgroup, then the orbits of $N$ in $\mathcal{X}$ are a partition and the orbits of $\{gNg^{-1}\}_{g \in G}$ form a partition base. One approach to the display of such projections is a 2-dimensional histogram using the ordering given by one of the metrics suggested in Chapter 7 of Diaconis [18].

## 2.3. **Projections for $\mathcal{X} = \mathbb{R}^p$: Euclidean data.**

Consider data vectors $x_1, x_2, \ldots, x_n \in \mathbb{R}^p$. For $\gamma$ in the $p$-dimensional unit sphere, the projection in direction $\gamma$ is just $\gamma \cdot x_1, \ldots, \gamma \cdot x_n$. This is the classical Radon transform, with $\mathcal{Y}$ consisting of the affine hyperplanes $y_\gamma^t = x \in \mathbb{R}^p : x \cdot \gamma = t$. For fixed $\gamma$ these partition the space $\mathbb{R}^p$ as $t$ varies, and the partitions vary as $\gamma$ varies. In statistical applications, a histogram is made of $\gamma \cdot x_i$ and one varies $\gamma$, trying to understand the structure of the $p$-dimensional data from the varying histograms. This leads to the classical version of projection pursuit considered in the introduction.

### 2.4. Projections when $\mathcal{X}$ is a finite set with $n$ elements, and $\mathcal{Y}$ is the class of $k$-element subsets.

In this example, if $k$ divides $n$, it is a non-trivial theorem of Baranyai that $\mathcal{Y}$ forms a projections base. Details and discussion may be found in Cameron [7]. This example occurs naturally when considering extensions of a given class of partitions. For example, consider the marginal projections $y_i^a$ in $\mathbb{Z}_2^k$ defined above. These sets all have cardinality $|y_i^a| = 2^{k-1}$. It is natural to consider the extension to projections based on the class of all subsets of cardinality $2^{k-1}$.

### 2.5. Uniqueness of Radon transforms:

We now consider the question: when is $f \to \bar{f}$ one to one? A convenient criteria involves the notion of a block design. Let $|\mathcal{X}| = n$. The class of sets $\mathcal{Y}$ is a block design with parameters $(n, c, k, l)$ provided

(4)                    $|y| = c$ for all $y \in \mathcal{Y}$,

(5)                    each $x \in \mathcal{X}$ is contained in $k$ subsets $y$,

(6)                    each pair $x \neq x^{'}$ is contained in $l$ subsets $y$.

Affine planes or $\mathbb{Z}_2^k$ and $k$ sets of an $n$ set are block designs. A great many other examples are discussed in the literature of combinatorial designs. In the statistics literature they are sometimes called balanced incomplete block designs. In the combinatorial literature they are often called 2-designs, or 2-$(n, c, l)$ designs. It is easy to see that the parameters $n, c, k, l$ satisfy

(7)                        $|\mathcal{Y}|c = nk,$

(8)                   $(n-1)l = k(c-1).$

Bailey [3], Dembroski [14] and Lander [38] are useful references for block designs. The following result is well known in the theory of designs. We first learned it from Bolker [4].

**Theorem 2.2.** *If $\mathcal{X}$ is a finite set and $\mathcal{Y}$ is a block design with $|\mathcal{Y}| > 1$, then the Radon transform $f \to \bar{f}$ is one to one, with an explicit inversion formula given by (12) below.*

*Proof.* For any $x$,

(9)          $\displaystyle \sum_{y: x \in y} \bar{f}(y) = kf(x) + l \sum_{\substack{s,s' \in \mathcal{X} \\ x \neq x'}} f(x')$

(10)                    $\displaystyle = (k-l)f(x) + l \sum_{x \in \mathcal{X}} f(x).$

If $\sum_{x \in \mathcal{X}} f(x) = 1$, this determines $f$ as

$$(11) \qquad f(x) = \frac{1}{k-l} \sum_{y:x \in y} \bar{f}(y) - \frac{l}{k-l}.$$

Observe that $k > l$ follows from the assumption that $|\mathcal{Y}| > 1$. When $\sum_{x \in \mathcal{X}} f(x)$ is not known, it can be recovered by summing both sides of (9) in $x$. This gives

$$\sum_{x \in \mathcal{X}} f(x) = \frac{c}{k-l+nl} \sum_{y \in \mathcal{Y}} \bar{f}(y)$$

and so the inversion formula

$$(12) \qquad f(x) = \frac{1}{k-l} \sum_{y:x \in y} \bar{f}(y) + \frac{lc}{(k-l)^2 + nl(k-l)} \sum_{y \in \mathcal{Y}} \bar{f}(y). \qquad \square$$

**Remarks.**

- It is not necessary that $\mathcal{Y}$ be a block design for $f \to \bar{f}$ to be one to one. For example, Kung [37] shows that the Radon transform is one to one when $\mathcal{Y}$ consists of the sets of rank $i$ in a matroid. Diaconis and Graham [17] give examples where the transform is one to one when $\mathcal{Y}$ consists of the nearest neighbors in a metric space. For example, when $\mathcal{X} = \mathbb{Z}_2^{2k}$ and $\mathcal{Y}$ consists of the balls of Hamming distance less than or equal to 1, the transform is one to one, and an explicit inversion theorem is known. When $\mathcal{X}$ is $S_n$, the symmetric group, and $\mathcal{Y}$ is unit balls in the Cayley metric, the transform is one to one if and only if $n$ is in $\{1, 2, 4, 5, 6, 8, 10, 12\}$. Further work on inversion formulas for functions on finite symmetric spaces is found in Velasquez [47] and for functions on the torus $\mathbb{Z}_n^k$ in Dedeo and Velasquez [13]. Fill [24] discusses invertibility when the Radon transform of $f$ at $x$ averages over a set of translates of $f(x)$ which has applications to directional data and time series.
- The transform can still be useful and interesting if it is not one to one. For example, the marginal projections in the example above do not capture all aspects of the data but are often the first things to be looked at. In $\mathbb{Z}_2^k$, if high enough marginal distributions are considered, the function $f$ can be completely recovered. In the symmetric group, the projections corresponding to all Young subgroups determine $f$ because they determine its Fourier transform. See Diaconis [18] for details.

## 3. Data analysis of syllable patterns in the works of Plato

This section presents a new analysis of data arising from syllable patterns in the works of Plato. The data are given in Table 1. It records, for each of 6 books, the pattern of long (-) and short ($\cup$) syllables among the last 5 syllables in each sentence. It is known that Plato wrote *Republic* early and *Laws* late. Plato also mentions that he changed his rhyming patterns over time. This led Brandwood to collect the data in Table 1.

The other books were written between these but it is not known in what order. The goal of the analysis is to try to order the books. Our approach will be to study the books one at a time, trying to find patterns.

Projection pursuit suggests looking at various partitions of the data, searching for structured partitions which are far from uniform. Using first and second order margins as partitions, a reasonably striking difference between *Republic* and *Laws* is observed. This suggests a simple, interpretable way of ordering the other books as *Republic, Timaeus, Sophist, Politicus, Philebus, Laws.*

This agrees with the standard ordering as discussed in Brandwood ([6], pg. xviii) and in Ahn et al. [1]. Other analyses of this data set are in Cox and Brandwood [11], Atkinson [2], Wishart and Leach [49], Boneva [5], and Charnomordic and Holmes [8]. [11] contains a history and explanation for the choice of data. The first three analyses all use statistical models. Boneva's analysis uses a form of scaling. None of the previous analyses seem to have picked up the simple, striking pattern in the data that projection pursuit leads to.

The analysis is presented below, in a somewhat discursive style, in the order it was actually performed: first looking at the *Republic*, then *Laws* and finally the other books. In the Appendix, we present a more automated and formal version.

### *3.1.* **Republic**

Table 2 shows the first order margins; e.g., the proportion of sentences with ∪ in position $i$, $1 \le i \le 5$.

Roughly, positions 1-4 are evenly divided between long and short. The last position is clearly different. Table 3 shows the second order margins.

A glance at Table 3 shows that the first order effects are all too visible in the second order margins. For example, the numbers in the first column (∪ ∪) are all "small" while the numbers in the last column are "large". One simple way of adjusting for the first order structure is to divide each number in Table 3 by the product of the marginal totals. For example, in the first row, .194 would be divided by (.465)(.472) (from Table 2) while .271 would be divided by (.465)(1 − .472). The results are shown in Table 4.

Most of the ratios are close to 1, so a product model is a reasonable first description. The projection pursuit approach suggests that a partition of the data (here a row) is "interesting" if the partition is far from uniform. By eye, looking at Table 4, positions $(1, 2), (2, 3), (3, 4), (4, 5)$ are far from being all 1. Observe that these positions are adjacent, as $(i, i + 1)$.

Next observe that each of the 4 designated rows has a common pattern: the first and last entries are small, the middle two entries are large. Going back to the definitions, this pattern arises from a negative association of adjacent syllables; in the *Republic*, adjacent syllables tend to alternate. The pattern in positions $(1,3)$ shows that this cannot be a complete description; after all, if the symbols alternate, the positions two apart should be positively associated, but $(1,3)$ displays negative association. Looking at the other rows of the table, we observe that the size goes big, small, small, big or its opposite, small, big, big, small. This is an artifact. Consider the first row of Table 4. It was formed from 4 proportions that sum to 1: $w, x, y, z$ say. The 4 adjusted entries are

$$\frac{w}{(w+x)(w+y)} \quad \frac{x}{(w+x)(x+z)} \quad \frac{y}{(y+z)(y+w)} \quad \frac{z}{(z+y)(z+x)}.$$

TABLE 2
*First order margins for* Republic

| Position | 1 | 2 | 3 | 4 | 5 |
|---|---|---|---|---|---|
| Proportion of ∪ | 0.465 | 0.471 | 0.466 | 0.511 | 0.362 |

TABLE 3
*Second order margins for* Republic

| Position | ∪∪ | ∪ - | - ∪ | - - |
|---|---|---|---|---|
| (1,2) | 0.194 | 0.271 | 0.277 | 0.258 |
| (1,3) | 0.208 | 0.257 | 0.258 | 0.277 |
| (1,4) | 0.238 | 0.227 | 0.272 | 0.263 |
| (1,5) | 0.177 | 0.288 | 0.185 | 0.350 |
| (2,3) | 0.209 | 0.262 | 0.257 | 0.272 |
| (2,4) | 0.241 | 0.230 | 0.269 | 0.260 |
| (2,5) | 0.162 | 0.309 | 0.200 | 0.329 |
| (3,4) | 0.211 | 0.255 | 0.299 | 0.235 |
| (3,5) | 0.170 | 0.296 | 0.192 | 0.342 |
| (4,5) | 0.167 | 0.343 | 0.195 | 0.295 |

It is easy to show that the first entry is less than 1 if and only if the second is larger than 1, if and only if the third is larger than one, if and only if the fourth is less than 1. This means that the first column in Table 4, together with the first order margins, determines the remaining entries. This artifact in no way reflects on the association pattern noted earlier– the most structured rows correspond to adjacent syllables, and adjacent syllables are negatively associated.

### 3.2. **Laws** *and a comparison with* **Republic.**

The first order margins for *Laws* are only slightly different from those in *Republic* (see Table 5).

The pattern is the same: overall, fewer than half ∪'s; the last position sharply smaller. The similarity between the first order margins in *Republic* and *Laws* suggests that second or higher order margins must be used to order the remaining books. The analog of the first column of Table 4 is given in Table 6.

The entries above are the proportion of sentences with ∪ ∪ in the $(i,j)$ position divided by the product of the marginal proportions.

TABLE 4
*djusted second order margins for* Republic

| Position | ∪ ∪ | ∪ - | - ∪ | - - |
|----------|------|------|------|------|
| (1,2) | 0.89 | 1.10 | 1.10 | 0.91 |
| (1,3) | 0.96 | 1.00 | 1.00 | 0.97 |
| (1,4) | 1.00 | 1.00 | 1.00 | 1.00 |
| (1,5) | 1.10 | 0.97 | 0.96 | 1.00 |
| (2,3) | 0.95 | 1.00 | 1.00 | 0.96 |
| (2,4) | 1.00 | 1.00 | 1.00 | 1.00 |
| (2,5) | 0.95 | 1.00 | 1.00 | 0.97 |
| (3,4) | 0.89 | 1.10 | 1.10 | 0.90 |
| (3,5) | 1.00 | 1.00 | 0.99 | 1.00 |
| (4,5) | 0.90 | 1.10 | 1.10 | 0.94 |

TABLE 5
*First order margins for* Laws

| Position | 1 | 2 | 3 | 4 | 5 |
|----------|------|------|------|------|------|
| Proportion of ∪ ∪ | 0.477 | 0.489 | 0.411 | 0.599 | 0.375 |

TABLE 6
*Adjusted second order margins for* Laws

| Positions | (1,2) | (1,3) | (1,4) | (1,5) | (2,3) |
|-----------|-------|-------|-------|-------|-------|
| Adjusted ∪ ∪ | 1.07 | 1.03 | 0.92 | 0.99 | 1.43 |
| Positions | (2,4) | (2,5) | (3,4) | (3,5) | (4,5) |
| Adjusted ∪ ∪ | 0.97 | 0.98 | 1.04 | 1.09 | 1.02 |

Again, pairwise adjacent positions are associated, all in the same way. Here, the association is positive, whereas for *Republic*, the association is negative. This is the striking pattern referred to above. It suggests a method of ranking the other books: compare the sign pattern or actual ratios of the adjusted second order margins of other books with *Republic* and *Laws*.

For definiteness, the sum of absolute deviations between second order margins over all 10 positions will be used. This is carried out data analytically in Sections 3.3–3.5.

### *3.3. Analysis for* Philebus *and* Politicus

These books are somewhat similar to each other. The first and second order margins for *Philebus* are given in Tables 7 and 8.

Note the difference in first order margins: between *Philebus* and *Republic* (or *Laws*) position 1 is high, as are positions 4 and 5. For second order margins, the adjacent patterns are all positively associated ((2,3) being truly extreme). Comparing Table 8 with Table 6, the association pattern matches *Laws* in direction, except in position (1,5). The relevant averages for *Politicus* are given in Tables 9 and 10.

The first order margins are, very roughly, like those in both *Republic* and *Laws*, but again the third position has a low proportion of short syllables. The second order margins have the same pattern as *Laws*. The same remarks made for the second order margins of *Philebus* apply.

Both *Philebus* and *Politicus* seem very similar to *Laws*. Which of these two is closest to *Laws*? One simple approach is to consider the sum of the absolute values of the differences between the entries of Tables 8 and 6 along with the differences between 10 and 6. The sum for *Laws* to *Philebus* is .64, while the sum for *Laws* to *Politicus* is .83. Thus a tentative ranking is: *Politicus*, *Philebus*, *Laws*.

TABLE 7
*First order margins for* Philebus

| Position | 1 | 2 | 3 | 4 | 5 |
|---|---|---|---|---|---|
| Proportion of ∪ | 0.522 | 0.464 | 0.398 | 0.594 | 0.465 |

TABLE 8
*Adjusted second order margins for* Philebus

| Positions | (1,2) | (1,3) | (1,4) | (1,5) | (2,3) |
|---|---|---|---|---|---|
| Adjusted ∪ ∪ | 1.11 | 1.03 | 0.85 | 1.11 | 1.48 |
| Positions | (2,4) | (2,5) | (3,4) | (3,5) | (4,5) |
| Adjusted ∪ ∪ | 0.92 | 0.85 | 1.02 | 0.95 | 1.01 |

TABLE 9
*First order margins for* Politicus

| Position | 1 | 2 | 3 | 4 | 5 |
|---|---|---|---|---|---|
| Proportion of ∪ | 0.477 | 0.457 | 0.348 | 0.524 | 0.469 |

TABLE 10
*Adjusted second order margins for* Politicus

| Positions | (1,2) | (1,3) | (1,4) | (1,5) | (2,3) |
|---|---|---|---|---|---|
| Adjusted ∪ ∪ | 1.17 | 1.10 | 0.96 | 1.01 | 1.26 |
| Positions | (2,4) | (2,5) | (3,4) | (3,5) | (4,5) |
| Adjusted ∪ ∪ | 0.86 | 0.90 | 1.05 | 1.10 | 1.13 |

### *3.4. Analysis for* Sophist *and* Timaeus

These books are quite similar to each other and, as we shall see, quite different from *Laws*, *Philebus* and *Politicus*.

The first order margins are quite different from the books examined previously. They are roughly consistent with all syllables being equally likely to be long or short. The first order pattern seems closest to *Politicus*. The second order associations are closer to 1 than in *Laws*, *Politicus* or *Philebus*. Adjacent positions are positively associated, except for (3,4). The direction of association matches *Laws* in only 6 of 10 positions. The sum of absolute deviations between the entries of Tables 6 and 12 is .87.

We now give the analysis for the final book.

A distinctive feature of the first order margins is the large proportion of short syllables in the third position. The adjusted second order margins are close to 1, so *Timaeus* seems closest to *Sophist*. Of the 4 adjacent positions, two show positive association and two show negative association. The direction of association matches *Laws* in 6 positions; the sum of absolute deviations between Tables 14 and 6 is .94. The distance between *Timaeus* and the *Republic* (Tables 14 and 4) is .6, so *Timaeus* seems closer to *Republic* than to *Laws* using this measure. Because of the decrease in the number of matches and the increase in the sum of absolute deviations, it seems reasonable to rank order the three as *Republic*, *Timaeus*, *Sophist*. This completes the discussion of this example. The Appendix contains an automated version.

## 4. Most projections are uniform

Graphical projection pursuit is a standard tool in data analysis. The classical survey of Huber [32], the survey of Posse [42] and the online documentation in the Xgobi and Ggobi packages contain extensive pointers to a large literature.

TABLE 11
*First order margins for* Sophist

| Position | 1 | 2 | 3 | 4 | 5 |
|---|---|---|---|---|---|
| Proportion of ∪ | 0.474 | 0.491 | 0.454 | 0.527 | 0.487 |

TABLE 12
*Adjusted second order margins for* Sophist

| Positions | (1,2) | (1,3) | (1,4) | (1,5) | (2,3)) |
|---|---|---|---|---|---|
| Adjusted ∪ ∪ | 1.07 | 1.03 | 1.01 | 0.93 | 1.07 |
| Positions | (2,4) | (2,5) | (3,4) | (3,5) | (4,5) |
| Adjusted ∪ ∪ | 0.88 | 1.01 | 0.97 | 0.98 | 1.10 |

TABLE 13
*First order margins for* Timaeus

| Position | 1 | 2 | 3 | 4 | 5 |
|---|---|---|---|---|---|
| Proportion of ∪ | 0.494 | 0.476 | 0.565 | 0.521 | 0.496 |

TABLE 14
*Adjusted second order margins for* Timaeus

| Positions | (1,2) | (1,3) | (1,4) | (1,5) | (2,3) |
|---|---|---|---|---|---|
| Adjusted ∪ ∪ | 0.98 | 1.02 | 0.97 | 1.04 | 0.92 |
| Positions | (2,4) | (2,5) | (3,4) | (3,5) | (4,5) |
| Adjusted ∪ ∪ | 0.94 | 0.97 | 0.96 | 0.97 | 1.06 |

The theorems of this section imply that for most data sets $f(x)$, most projections $\bar{f}(y)$ are about the same: close to uniform. This necessitates projection pursuit – choosing projections that are far from uniformly distributed – to determine what is special about a particular $f$. This gives an independent rationale for Huber's suggestion that Euclidean projections are interesting if they are far from uniform in the sense of having minimum entropy (of course, the uniform distribution on a finite set has maximum entropy).

**Theorem 4.1.** *Let $\mathcal{X}$ be a finite set with $n$ elements. Let $\mathcal{Y}$ be a block design with block size $c$ (so $|y| = c$ for $y \in \mathcal{Y}$). Let $f : \mathcal{X} \to \mathbb{R}$ be any function and let $\mu(f) = \sum_{x \in \mathcal{X}} f(x)$. Let $y$ be chosen uniformly in $\mathcal{Y}$. Then*

$$\text{(13)} \qquad \mathbf{E}\bar{f}(y) \;=\; \frac{c}{n}\mu(f),$$

$$\text{(14)} \qquad \text{var}\,\bar{f}(y) \;=\; \frac{c}{n}(1 - \frac{(c-1)}{(n-1)})\mu(f - \frac{\mu(f)}{n})^2.$$

*Proof.* (13) follows from computing

$$\mathbf{E}\bar{f}(y) = \frac{1}{|\mathcal{Y}|}\sum_y \bar{f}(y) = \frac{1}{|\mathcal{Y}|}\sum_x f(x)|y : x \in y| = \frac{c}{|\mathcal{Y}|}\mu(f).$$

For (14), assume without loss of generality, that $\mu(f) = 0$. Then

$$
\begin{aligned}
\text{var}(\bar{f}(y)) &= \frac{1}{|\mathcal{Y}|}\sum_y \bar{f}(y)^2 = \frac{1}{|\mathcal{Y}|}\sum_y \left( \sum_{x \in y} f(x)(f(x) + \sum_{\substack{x \neq x' \\ x, x' \in y}} f(x')) \right) \\
&= \frac{k-l}{|\mathcal{Y}|}\mu(f^2).
\end{aligned}
$$

From (7) and (8), $\frac{k-l}{|\mathcal{Y}|} = \frac{c(n-c)}{n(n-1)}$, giving the result. $\qquad\square$

**Example 4.2.** When $\mathcal{Y}$ is the $j$ sets of an $n$ set, $|\mathcal{Y}| = \binom{n}{j}, c = j$, and the result reduces to the usual mean and variance for a sample without replacement.

**Example 4.3.** Let $\mathcal{X} = \mathbb{Z}_2^k$ and $\mathcal{Y}$ be the j-dimensional affine planes. Then $n = 2^k$ and $c = 2^{k-j}$. If $\mu(f) = 1$, the result becomes

$$\mathbf{E}(\bar{f}(y)) = \frac{1}{2^j}, \qquad \text{var}(\bar{f}(y)) = \frac{1}{2^j}(1 - \frac{2^{k-j}-1}{2^k-1})\mu(f - \frac{1}{2^k})^2.$$

For future use, observe that the cardinality of $\mathcal{Y}$ in this case is

$$\frac{2^j(2^k-1)(2^k-2)\cdots(2^k-2^{j-1})}{(2^j-1)\cdots(2^j-2^{j-1})}.$$

Returning to the situation in Theorem 2.2, Chebychev's inequality implies:

**Corollary 4.4.** *With notation as in Theorem 2.2, the proportion of $y \in \mathcal{Y}$ such that*

$$|\bar{f}(y) - \frac{c}{n}\mu(f)| > \epsilon$$

*is smaller than*

$$\frac{1}{\epsilon^2}\frac{c}{n}(1 - \frac{c-1}{n-1})\mu(f - \frac{\mu(f)}{n})^2.$$

**Remarks.** The corollary implies that for functions $f$ which are "not too wild" in the sense that $\mu(f - \frac{\mu(f)}{n})^2$ is small, most transforms $\bar{f}(y)$ are uninformative in the sense of being close to their mean value. As an example, take $\mathcal{X} = \mathbb{Z}_2^5$ and $f$ the function defined by the first column of Table 1. Then $\mu(f - \frac{1}{32})^2 = .0021$. If $\mathcal{Y}$ is taken as the set of all affine hyperplanes, the corollary gives that 95% of the transforms have $|\bar{f}(y) - \frac{1}{2}| < .04$.

The next theorem says that for most probabilities $f$, $\mu(f - \frac{1}{n})^2$ is small (about $\frac{1}{n}$).

**Theorem 4.5.** *Let $(U_1, U_2, \ldots, U_n)$ be chosen uniformly on the $n$ simplex. For large $n$, the random variable*

$$\frac{n^{3/2}}{2} \left( \sum_{i=1}^{n} (U_i - \frac{1}{n})^2 - \frac{1}{n} \right)$$

*has an approximate standard normal distribution.*

*Proof.* The argument uses the representation of a uniform distribution by means of exponential variables. Let $X_1, X_2, \ldots, X_n$ be independent standard exponential variables with density $e^{-x}$ on $[0, \infty)$. Let

$$S_1 = \sum_{i=1}^{n} X_i, \qquad S_2 = \sum_{i=1}^{n} X_i^2.$$

For large $n$, the random vector

$$\binom{Z_1}{Z_2} = \frac{1}{\sqrt{n}} \binom{S_1 - n}{S_2 - 2n}$$

has an approximate bivariate normal distribution with mean vector zero and co-variance matrix $\left( \begin{smallmatrix} 1 & 4 \\ 4 & 20 \end{smallmatrix} \right)$. To check the covariance matrix, note that $\text{var}(\frac{S_1 - n}{\sqrt{n}}) = \text{var}(X_1) = 1$, $\text{var}(\frac{S_2 - n}{\sqrt{n}}) = \text{var}(X_1^2) = 20$ and

$$\begin{aligned} \frac{1}{n}\mathbf{E}((S_1 - n)(S_2 - 2n)) &= \mathbf{E}\left( (X_1 - 1)(X_1^2 - 2) \right) \\ &= \mathbf{E}(X_1^3) - \mathbf{E}(X_1^2) - 2\mathbf{E}(X_1) + 2 = 4. \end{aligned}$$

Represent a uniform vector on the $n$ simplex as $U_i = \frac{Z_i}{S_1}$. Then

$$\sum_{i=1}^{n} (U_i - \frac{1}{n})^2 = \frac{1}{S_1^2} \sum_{i=1}^{n} X_i^2 - \frac{1}{n} = \frac{1}{S_1^2} \sum_{i=1}^{n} (X_i^2 - 2) + \frac{2n}{S_1^2} - \frac{1}{n}.$$

Now $S_1 = n(1 + \frac{Z_1}{\sqrt{n}})$ with $Z_1 = \frac{S_1 - n}{\sqrt{n}}$. Thus

$$S_1^2 = n^2(1 + \frac{2}{\sqrt{n}}Z_1 + \frac{Z_1^2}{n}).$$

Using the standard $O_p$ notation (see Pratt [43]),

$$\frac{1}{S_1^2} = \frac{1}{n^2} - \frac{2Z_1}{n^{3/2}} + O_p(\frac{1}{n^3}).$$

Thus,

$$\frac{1}{S_1^2}\sum_{i=1}^n (X_i^2 - 2) = \frac{1}{n^{3/2}}\frac{1}{\sqrt{n}}\sum_{i=1}^n (X_i^2 - 2) + O_p(\frac{1}{n^2}),$$

$$\frac{2n}{S_1^2} = \frac{2}{n} - \frac{4Z_1}{n^{3/2}} + O_p(\frac{1}{n^2}).$$

The bivariate limiting normality of $\binom{Z_1}{Z_2}$ implies that $Z_2 - 4Z_1$ has an approximate normal distribution with mean 0 and variance

$$\text{var}(Z_2) + 16 \text{ var}(Z_1) - 8 \text{ covar}(Z_1, Z_2) = 4. \qquad \square$$

Corollary 4.4 and Theorem 4.5 imply that for most probabilities $f$, most transforms $\bar{f}(y)$ are close to uniform.

The final result of this section deals with the entire projection $\bar{f}(y)_{y \in p}$ where $p$ is a partition of $\mathcal{X}$ into blocks in $\mathcal{Y}$. Let $\mathcal{X}$ be a finite set. Let $\mathcal{Y}$ be a block design on $\mathcal{X}$ with parameters $(n, c, k, l)$. Suppose that $\mathcal{Y}$ is also a projection base for $\mathcal{X}$ with $p_1, p_2, \ldots, p_j$ being a partition of $\mathcal{Y}$, and each $p_i$ being a partition of $\mathcal{X}$. Of course, $j = \frac{|\mathcal{Y}|c}{n}$. The next theorem implies that for most functions, the projection onto a randomly chosen partition is uniformly close to $\frac{c}{n}$.

**Theorem 4.6.** *Let $\mathcal{Y}$ be a block design on $\mathcal{X}$ with parameters $(n, c, k, l)$. Suppose that $\mathcal{Y}$ is a projection base. Let $f$ be a fixed probability on $\mathcal{X}$. Let the partition $p$ be chosen uniformly at random over all partitions $p_i$ of $\mathcal{X}$, where $p_i \subset \mathcal{Y}$. For $\epsilon > 0$,*

$$(15) \qquad \sum_{y \in p} |\bar{f}(y) - \frac{c}{n}| \le \epsilon.$$

*with probability at least*

$$1 - \frac{1}{\epsilon}\left(\frac{n(n-c)}{c(n+1)}\mu(f - \frac{1}{n})^2\right)^{\frac{1}{2}}.$$

*Proof.* The probability model for choosing a random partition is based on a fixed enumeration $p_1, p_2, \ldots, p_j$ of the partitions that make up $\mathcal{Y}$. Each partition is assumed to be taken in a fixed order $p_i = \{(y_i^1, \ldots, y_i^{n/c})\}$. The random variable $S(p) = \sum_{y \in p}|\bar{f}(y) - \frac{c}{n}|$ is invariant under permuting the $y \in p$ among themselves. Thus a random variable with the same distribution of $S(p)$ but exchangeable $\bar{f}(y)_{y \in p}$ exists. For this realization, $\mathbf{E}(\sum_{y \in p}|\bar{f}(y) - \frac{c}{n}|) = \frac{n}{c}\mathbf{E}|\bar{f}(y^*) - \frac{c}{n}|$ with $y^*$ chosen uniformly in $\mathcal{Y}$. Using Cauchy-Schwartz and Theorem 4.1, the expectation is bounded above by

$$\frac{n}{c}\sqrt{\frac{c}{n}\left(1 - \frac{c-1}{n-1}\right)\mu(f - \frac{1}{n})^2}.$$

Theorem 4.6 follows from this bound and Markov's inequality applied to the original random variable. $\qquad \square$

**Remarks.** From Theorem 4.5, $\mu(f - \frac{1}{n})^2 \doteq \frac{1}{n}$ for most functions $f$. For such $f$, the theorem implies that for large block size $c$, most partitions are close to uniform in variation distance. This may be contrasted with Theorems 4.1 and 4.5 which imply that the components $\bar{f}(y)$ of most projections are close to $\frac{c}{n}$. When $c$ is small, there are many terms in the sum (15). As an example, consider the 2-sets of an $n$ set where $n = 2j$. Let $p$ be a random partition into 2-element sets. Let $f$ be chosen at random from the $n$ simplex and $p$ any fixed partition into two element sets. It is straightforward to show that with probability tending to 1 as $n$ tends to infinity,

$$\sum_{y \in p} |\bar{f}(y) - \frac{2}{n}| \to 8e^{-2}.$$

The analogous result holds with the same assumptions when $p$ is any fixed partition of fixed size $c$. Similarly, it is natural to ask for a central limit theorem in connection with Theorems 4.1 and 4.5. For $j$ sets of an $n$ set, such a theorem is available from the usual results on sampling without replacement from a finite population. Most likely, there is a similar set of results for block designs with $|\mathcal{Y}|$ and $c$ large. See Stein [45] for results for designs arising from subgroups of a finite group.

## 5. Least uniform partitions

The results of Section 4 imply that, under suitable conditions, for most functions the projection along most partitions is close to uniform. This suggests that the special properties of particular functions are only seen in partitions that are far from uniform. In this section, properties of least uniform partitions are examined. Theorem 5.1 shows that for most functions, even the least uniform partitions will be close to uniform if the the number of sets in $\mathcal{Y}$ is small in the sense that $\log |\mathcal{Y}|$ is small compared both to $n$ and the block size $c$. This is true, in particular, for affine hyperplanes in $\mathbb{Z}_2^k$.

**Theorem 5.1.** *Let $\mathcal{X}$ be a set of $n$ elements. Let $\mathcal{Y}$ be a class of subsets in $\mathcal{X}$ of fixed cardinality $c$. Suppose that $p_1, \ldots, p_j$ is a partition of $\mathcal{Y}$ into partitions of $\mathcal{X}$. Let $f$ be chosen at random in the $n$ simplex. Let $p^*$ be the partition in $p_i$ that maximizes $\sum_{y \in p} |\bar{f}(y) - \frac{c}{n}|$. For any $\epsilon > 0$,*

$$\sum_{y \in p^*} |\bar{f}(y) - \frac{c}{n}| < \epsilon,$$

*except for a set of $f$'s of probability smaller than*

$$(|\mathcal{Y}| + 1)\beta$$

*with $\beta$ equal to 1 minus*

$$(16) \qquad \frac{1}{\beta(c, n)} \int_{\frac{c}{n}(1-\epsilon)}^{\frac{c}{n}(1+\epsilon)} x^{c-1}(1-x)^{n-c-1} dx,$$

*where $\beta(c, n)$ denotes the beta function.*

*Proof.* Represent the $i^{\text{th}}$ component of a randomly chosen $f$ as $\frac{X_i}{S}$ where $X_i$ are independent standard exponentials and $S = \sum_{i=1}^n X_i$. Let $y^*$ be the set in $\mathcal{Y}$ with the largest value of $\frac{c}{n}(1-\epsilon)$. The argument begins by bounding the probability that

$$|\bar{f}(y^*) - \frac{c}{n}| < \epsilon \frac{c}{n}.$$

To begin with,

$$\mathbf{P}\left(\bar{f}(y^*) < \frac{c}{n}(1-\epsilon)\right) \leq \mathbf{P}\left(\frac{X_1 + \cdots + X_c}{S} < \frac{c}{n}(1-\epsilon)\right).$$

Further,

$$\mathbf{P}\left(\bar{f}(y^*) > \frac{c}{n}(1+\epsilon)\right) \leq \sum_{y \in \mathcal{Y}} \mathbf{P}\left(\bar{f}(y) > \frac{c}{n}(1+\epsilon)\right)$$

$$= |\mathcal{Y}| \mathbf{P}\left(\frac{X_1 + \cdots + X_c}{S} > \frac{c}{n}(1+\epsilon)\right).$$

Next, let $y_*$ denote the set in $\mathcal{Y}$ with the smallest value of $\bar{f}(y)$. To bound the probability that $|\bar{f}(y_*) - \frac{c}{n}| < \epsilon\frac{c}{n}$, observe that $\bar{f}(y_*) = 1 - \bar{f}(y^{**})$ with $y^{**}$ the union of sets in a partition omitting the one element that maximizes $\bar{f}$. Thus,

$$\mathbf{P}\left(\bar{f}(y_*) < \frac{c}{n}(1-\epsilon)\right) = \mathbf{P}\left(\bar{f}(y^{**}) > 1 - \frac{c}{n}(1-\epsilon)\right)$$

$$\leq |\mathcal{Y}| \mathbf{P}\left(\frac{X_1 + \cdots + X_{n-c}}{S} > 1 - \frac{c}{n}(1-\epsilon)\right)$$

$$= |\mathcal{Y}| \mathbf{P}\left(\frac{X_1 + \cdots + X_c}{S} < \frac{c}{n}(1-\epsilon)\right).$$

Further,

$$\mathbf{P}\left(\bar{f}(y_*) > \frac{c}{n}(1+\epsilon)\right) = \mathbf{P}\left(\bar{f}(y^{**}) > 1 + \frac{c}{n}(1-\epsilon)\right)$$

$$\leq \mathbf{P}\left(\frac{X_1 + \cdots + X_{n-c}}{S} < 1 - \frac{c}{n}(1+\epsilon)\right)$$

$$= \mathbf{P}\left(\frac{X_1 + \cdots + X_c}{S} > \frac{c}{n}(1+\epsilon)\right).$$

Summing the four bounds thus obtained we see that both

(17) $$\qquad |\bar{f}(y_*) - \frac{c}{n}| < \epsilon\frac{c}{n}, \qquad\qquad |\bar{f}(y^*) - \frac{c}{n}| < \epsilon\frac{c}{n}$$

except for a set of $f$'s of probability smaller than $(|\mathcal{Y}|+1)\beta$ as defined by (16). Now (17) implies that $|\bar{f}(y) - \frac{c}{n}| < \epsilon\frac{c}{n}$ for all $y \in \mathcal{Y}$. Summing this last inequality over the partition $p^*$ completes the proof of the theorem. $\qquad\square$

**Remarks.** The beta integral that appears in the bound is straightforward to approximate numerically. A raft of techniques and approximations appear in the first chapter of Pearson [40]. For example, consider cases where $\frac{c}{n} = \frac{1}{2}$. Then, using the Peizer-Pratt approximation given in Pearson [40], and Mills' ratio, the $\beta$ in (16) is approximately

$$\frac{2}{\sqrt{2\pi}} \frac{e^{-\frac{x^2}{2}}}{1+x} \qquad \text{with } x = \sqrt{2c \log \frac{1}{4(\frac{1}{2}-\epsilon)(\frac{1}{2}+\epsilon)}}.$$

For this to be small when multiplied by $|\mathcal{Y}| + 1$, it clearly suffices that $\log|\mathcal{Y}|$ be small compared to $c$. This is the case for the affine subspaces of dimension $j$ in $\mathbb{Z}_2^k$ if $j$ is bounded and $k$ is large.

As a numerical example, consider the affine hyperplanes in $\mathbb{Z}_2^{10}$. Then $|\mathcal{Y}| + 1 = 2049$, $c = 512$, $n = 1024$. Taking $\epsilon = 0.1$, $(|\mathcal{Y}| + 1)/\beta \doteq 2.595 \times 10^{-7}$.

The next theorem shows that when there are many sets in $\mathcal{Y}$, the least uniform projection is typically far from uniform. The theorem deals with $n$ sets in a set of cardinality $2n$. The variation distance of a typical probability projected along the least uniform half split is shown to be about 0.3. This may be compared with Theorems 4.5 and 4.6 which show that for a typical probability $f$ on $2n$ points, $|\bar{f}(y) - \frac{1}{2}|$ is close to zero for most sets $y$ of cardinality $n$.

**Theorem 5.2.** *Let $f$ be chosen at random on the $2n$ simplex. Let $S^-$ be the sum of the $n$ smallest $f(x)$. Then for large $n$, the random variable*

$$\sqrt{2n}\left(S^- - (\frac{1}{2} - \frac{\log 2}{2})\right)$$

*has an approximate normal distribution with mean 0 and variance $\frac{3}{2} - 2\log 2$.*

*Proof.* Represent a randomly chosen $f$ as $\frac{X_i}{S}$ where $X_i$ are independent standard exponential random variables and $S = \sum_{i=1}^{2n} X_i$. Denote the order statistics by round brackets:

$$X_{(1)} \leq X_{(2)} \leq \cdots \leq X_{(n)}.$$

Let $L_1 = X_{(1)}, L_2 = X_{(2)} - X_{(1)}, \ldots, L_{2n} = X_{(2n)} - X_{(2n-1)}$. Then the $L_i$ are independent, and $L_{i+1}$ has the distribution of a standard exponential times $\frac{1}{(2n-i)}$ – see Feller ([23] Section III.3). With this notation,

$$S = \sum_{i=1}^{2n} X_i = \sum_{i=0}^{2n-1} (2n-i)L_{i+1}, \tag{18}$$

$$S^- = \frac{1}{S}\sum_{i=1}^{n} X_{(i)} = \frac{1}{S}\sum_{i=0}^{n-1}(n-i)L_{i+1}. \tag{19}$$

The proof is completed by approximating the sums in this representation of $S$ and $S^-$. Let $\mu_i = \frac{n-i}{2n-i}$, so $(n-i)L_{i+1}$ has the same distribution as $\mu_i$ times a standard exponential. Let

$$\begin{aligned}
\sigma^2 = 2\sum_{i=0}^{n-1}\mu_i^2 &= 2\sum_{i=0}^{n-1}(1 - \frac{2n}{2n-i} + \frac{n^2}{(2n-i)^2}) \\
&= 2\left(n - (2n\log 2 + O(1)) + \frac{3}{2} + \frac{n}{2} + O(1)\right) \\
&= 2n\left(\frac{3}{2} - 2\log 2\right) + O(1).
\end{aligned}$$

Now, let $Z_1 = \frac{S - 2n}{\sqrt{2n}}$ and $Z_2 = \frac{(\sum_{i=1}^{n} X_{(i)} - \mu_i)}{\sqrt{2n}}$. The vector $(Z_1, Z_2)$ has a limiting bivariate normal distribution, with mean $(0,0)$ and covariance matrix $\begin{pmatrix} \sigma_1^2 & \rho \\ \rho & \sigma_2^2 \end{pmatrix}$ with $\sigma_1^2 = 2, \sigma_2^2 = \frac{3}{2} - 2\log 2$, and $\rho = \frac{1}{2}(1 - \log 2)$. To check the value of $\rho$, observe that the covariance of $Z_1$ and $Z_2$ is $\frac{1}{2n}$ times

$$\begin{aligned}
\sum_{i=0}^{n}\mathbf{E}\left(\left[(2n-i)L_{i+1} - 1\right]\left[(n-i)L_{i+1} - \frac{1}{2n-i}\right]\right) &= \sum_{i=0}^{n}\frac{n-i}{2n-i} \\
&= n - n\log 2 + O(1).
\end{aligned}$$

Using the standard $O_p$ calculus,

$$\frac{1}{S} = \frac{1}{2n}\frac{1}{(1+\frac{Z_1}{\sqrt{2n}})} = \frac{1}{2n}\left(1 - \frac{Z_1}{\sqrt{2n}}\right) + O_p\left(\frac{1}{n^2}\right).$$

In particular,

$$\frac{1}{S} = \frac{1}{2n} + O_p\left(\frac{1}{n^{\frac{3}{2}}}\right).$$

The representation (19) for $S^-$ can be rewritten as

$$S^- = \sqrt{2n}\,\frac{Z_2}{X} + \frac{\mu}{S} = \frac{Z_2}{\sqrt{2n}} + \frac{1-\log 2}{2}\left(1 - \frac{Z_1}{\sqrt{2n}}\right) + O_p\left(\frac{1}{n}\right).$$

It follows that $\sqrt{2n}(S^- - \frac{1-\log 2}{2})$ has the same limiting distribution as $Z_2 - \frac{(1-\log 2)}{2}Z_1$. This is normal with mean 0 and variance

$$\left(\frac{3}{2} - 2\log 2\right) + 2\left(\frac{1-\log 2}{2}\right) - 2\left(\frac{1-\log 2}{2}\right) = \frac{3}{2} - 2\log 2. \qquad \square$$

**Corollary 5.3.** *Let $f$ be chosen at random on the $2n$ simplex. Let $(y, y^c)$ be a partition of $\mathcal{X}$ into an $n$ set and its complement which maximizes the value of*

$$|\bar{f}(y) - \frac{1}{2}| + |\bar{f}(y^c) - \frac{1}{2}|.$$

*Then, as $n$ tends to infinity, the maximum discrepancy tends to $\log 2 \doteq .301$ with probability tending to 1.*

*Proof.* For almost all $f$, the maximum is taken on uniquely at the partition $S^-$, $(S^-)^c$ as defined in Theorem 5.2. The maximum discrepancy equals

$$2|S^- - \frac{1}{2}|,$$

and the result follows from Theorem 5.2. $\qquad \square$

**Remark.** The proof of Theorem 5.2 and its corollary can easily be extended to cover the $j$ sets of an $n$ set. The argument shows that for most probabilities $f$, the variation distance between the least uniform projection and the uniform distribution is bounded away from zero if $j$ is an appreciable fraction of $n$.

For the final theorem, a different method of choosing a random probability is introduced. Let $\mathcal{X}$ be a set of cardinality $2n$. Fix an integer $b$. Drop $b$ balls into $2n$ boxes, and let $f(x)$ be the proportion of balls in the box labeled $x$. Let $\mathcal{Y}$ be the subsets of $\mathcal{X}$ with cardinality $n$. Clearly, if $b$ is large with respect to $n$, $f(x)$ is approximately $\frac{1}{2n}$ and so for any $y \in \mathcal{Y}$, $\bar{f}(y) \doteq \frac{1}{2}$, even for the $y_*$ minimizing $\bar{f}(y)$. At the other extreme, if $b$ is small with respect to $n$, $\bar{f}(y_*)$ will be close to zero. For example, if $b = n$, $\bar{f}(y_*) = 0$. It will follow from Theorem 5.4 that $\bar{f}(y_*)$ is approximately zero for $v \le 2n\log 2$.

This model for generating a random probability gives insight into the following problem. If data is generated from a structureless model, random fluctuations may produce structure that is picked up by a rich enough data analytic procedure. As $b$ varies in the above model, the random probability converges to a uniform distribution. The following theorem gives an indication of how large $b$ must be for all projections to be close to uniform. Some required notation: For $\lambda < 0$, let $p_\lambda(j) = \frac{e^{-\lambda}\lambda^j}{j!}$ denote the Poisson density. Let $P_\lambda(j) = \sum_{i=0}^{j} p_\lambda(i)$. Let $m$ be the largest integer with $P_\lambda(m) \leq \frac{1}{2}$, $P_\lambda(m+1) > \frac{1}{2}$. Define $\theta = \theta(\lambda)$ by

$$P_\lambda(m) + \theta p_\lambda(m+1) = \frac{1}{2}, \qquad \text{so } 0 \leq \theta < 1.$$

Table 15

| $\lambda$ | 1 | 2 | 3 | 4 | 5 | 6 | 7 | 8 | 9 | 10 |
|---|---|---|---|---|---|---|---|---|---|---|
| $\frac{2e^{-\lambda}\lambda^m}{m!}$ | 0.74 | 0.54 | 0.44 | 0.40 | 0.36 | 0.32 | 0.30 | 0.28 | 0.26 | 0.24 |

When $\lambda$ is an integer, Ramanujan showed that

$$\theta = \frac{1}{3} + O(\frac{1}{\lambda}) \text{ as } \lambda \to \infty.$$

See Cheng [9] for references and extensions of Ramanujan's results.

**Theorem 5.4.** *Suppose that $n$ and $b$ tend to infinity in such a way that $\frac{b}{2n} \to \lambda$. Let $y_*$ be the $n$ set with smallest value of $\bar{f}(y_*)$. Then*

$$|\bar{f}(y_*) - \frac{1}{2}| + |\bar{f}(y_*^c) - \frac{1}{2}| = \frac{2e^{-\lambda}\lambda^m}{m!}\left(1 + \theta(\frac{\lambda}{m+1} - 1)\right) + o_p(1).$$

**Remarks.** For $\lambda \leq \log 2$ and $m = 0$ , the variation distance can be shown to tend to one. For large $\lambda$, $\frac{e^{-\lambda}\lambda^m}{m!}$ is roughly $\frac{1}{\sqrt{2\pi\lambda}}$; thus for large $\lambda$, the variation distance tends to zero like $\frac{1}{\sqrt{\lambda}}$. This is not very rapid as Table 15 shows. (Note that for integer $\lambda$, $m + 1 = \lambda$, so the asymptotic value of the variation distance is $2\frac{e^{-\lambda}\lambda^m}{m!}$.)

*Proof.* The argument will only be sketched. For $b$ and $n$ large, the number of balls in the $i^{\text{th}}$ box has a limiting Poisson distribution with parameter $\lambda$, and different boxes can be treated as independent. The arguments in Diaconis and Freedman ([15], Section 3) can be used to justify this step.

Thus let $X_1, X_2, \ldots$ be independent Poisson variables with mean $\lambda$. With probability 1, eventually the median of $X_1, X_2, \ldots, X_{2n}$ is $m + 1$ and the proportion of $X_i$, $1 \leq i \leq 2n$ equal to $j$ is $p_\lambda(j) + o(1)$ uniformly for $0 \leq j \leq m + 1$. Let $S^-$ be the sum of the $n$ smallest $X_i$, $1 \leq i \leq 2n$. It follows that $\frac{S^-}{2n}$ equals

$$0p_\lambda(0) + p_\lambda(1) + \cdots + mp_\lambda(m) + \theta(m+1)p_\lambda(m+1) + o(1).$$

This sum equals

$$\frac{\lambda}{2} - \frac{e^{-\lambda}\lambda^m}{m!}\left(1 + \theta(\frac{\lambda}{m+1} - 1)\right) + o(1).$$

The identity asserted in the theorem follows from noting that $\bar{f}(y_*)$ is the limiting value of $\frac{S^-}{2\lambda n}$. □

## Appendix: Automating the analysis

In Section 2, we used the adjusted second order margins in a graphical, data analytic fashion to seriate the books of Plato. For some purposes, it may be desirable to have a more formal ranking procedure. We carry this out in Section A.1. The procedure is based on a collection of metrics between probabilities. These are explained in Section A.2. Finally, in Section A.3, we carry out a fully automated analysis of the Plato data based on all affine projections, not just first and second order statistics. We conclude that most methods agree, and suggest that the structures described in Section 3 are robustly embedded in the Plato data. In this section, we have added a seventh book, Criticus, to the analysis.

### A.1. *A metric approach*

In our data analysis, the adjusted second order statistics emerged as an informative summary of the rhyming patterns in Plato's *Republic*. As explained in Section 2, this is a vector of ten numbers (one for each pair of the last five syllables, i.e. $\binom{5}{2} = 10$). For the moment, call this vector $p^R = (p_1^R, \ldots, p_{10}^R)$ with "$R$" denoting *Republic*. A similar ten-vector can be computed for each of the other books. We may then use the distance between these vectors and $p^R$ to order the books. Books closest to $p^R$ are ranked earlier. We also compute a ranking based on the distance to $p^L$, the adjusted second order statistics for Plato's *Laws*. These two rankings generally agree, and agree with the conclusions of Section 3.

To proceed, we need to choose a distance between vectors. We have examined three standard distances between probability vectors: the Hellinger distance $H$, the Total Variation distance $TV$, and the Vasserstein distance $V$. These are explained more carefully in Section A.2. The rankings are given in Table 16: $R$ denotes *Republic*, $L$ denotes *Laws*, · denotes row variable.

Almost the same seriation is obtained when any of the three metrics are used to compute distances between *Republic* and the other books. Similarly, almost the same seriation is obtained when any of the three metrics are used to compute distances between *Laws* and the other books. Most clearly, *Politicus* is closest to *Laws* and furthest from the *Republic*. *Timaeus* and *Sophist*, as a pair, are closest to *Republic* and furthest from *Laws*. However, *Sophist* is both closer to *Laws* and to *Republic* than *Timaeus*. From these calculations, aside from *Politicus*, *Philebus* is closest to *Laws* and furthest from *Republic*. This is then followed by *Criticus*. All of this points to the ordering: *Republic*, {*Sophist*, *Timaeus*}, *Criticus*, *Philebus*, {*Politicus*, *Laws*}.

This ordering is consistent with the ordering produced data analytically in Section 3 and with the ordering based on the exponential model of Cox and Brandwood [11]. In Ahn et al. [1], a total of ten books were used for analysis. They found "roughly three clusters" (618): {*Tim.*, *Soph*, *Crit.*, *Pol.* \* } { *Laws*, *Phil.* }, { *Rep*, \*,\* }. Here \* denotes a book not analyzed in our work. Their final ordering based on a cluster analysis using the Euclidean metric is *Republic*, *Timaeusus*, *Criticus*, *Sophist*, *Politicus*, *Philebus*, *Laws*.

TABLE 16
*Ranking of book in row based on distance in column*

| Book | $d_H(R, \cdot)$ | $d_{TV}(R, \cdot)$ | $d_V(R, \cdot)$ | $d_H(L, \cdot)$ | $d_{TV}(L, \cdot)$ | $d_V(L, \cdot)$ |
|------|------|------|------|------|------|------|
| Tim. | 2 | 2 | 2 | 5 | 5 | 5 |
| Soph. | 1 | 1 | 1 | 4 | 4 | 4 |
| Pol. | 6 | 5 | 6 | 1 | 1 | 1 |
| Crit. | 3 | 3 | 3 | 3 | 3 | 3 |
| Phil. | 4 | 4 | 4 | 2 | 2 | 2 |
| Laws. | 5 | 6 | 5 | 0 | 0 | 0 |

### *A.2. Some metrics*

Let $p = (p_1, \ldots p_n)$, $q = (q_1, \ldots q_n)$ be probability vectors. Thus $p_1 \geq 0$ and $p_1 + \ldots + p_n = 1$, and the same holds for $q$. Three widely used metrics are :

| | | |
|---|---|---|
| Total Variation: | $d_{TV}(p, q) = \frac{1}{2} \sum_i |p_i - q_i|$. | |
| Hellinger: | $d_H(p, q) = \sum_i (\sqrt{p_i} - \sqrt{q_i})^2$. | |
| Vasserstein: | $d_V(p, q) = \min_{X,Y} \mathbf{E}(d(X, Y))$. | |

where the minimum is over all joint distributions of $X$ and $Y$ with marginals $p$ and $q$.

These metrics, their strengths, weaknesses and relations are discussed in Dudley [21], Villani [48] and Diaconis et al. [20].

In Section A.1, we used these metrics between vectors of positive entries which did not necessarily have sum one. This was done by forming $\bar{p} = \sum_i p_i$, $\bar{q} = \sum_i q_i$, $\tilde{p} = \frac{p_i}{\bar{p}}$, $\tilde{q} = \frac{q_i}{\bar{q}}$. We used the distance between $\tilde{p}$ and $\tilde{q}$ and added a penalty term to account for differences in mass between the profiles $p$ and $q$. For total variation, the penalty was $|\bar{p} - \bar{q}|$. We computed and compared two penalty terms for Hellinger: both $|\bar{p} - \bar{q}|$ and $(\sqrt{\bar{p}} - \sqrt{\bar{q}})^2$.

Thus, the distances between the ten-vector of adjusted second order margins of *Republic* and the other books, using Vasserstein is given in Table 17.

For completeness, we note that the Vasserstein metric requires an underlying distance on a probability space; in our case, this amounts to an underlying distance between the ten entries in each table. We take these entries to be binary 5-tuples containing two ones. We use the distance between two of these as the minimum number of pairwise adjacent switches required to bring one to the other. Thus the distance between 11000 and 00011 is 6. Further background can be found in Diaconis et al. [20] or Thompson [46]. With this choice specified, the minimization problem is equivalent to the Monge-Kantorovich Transshipment problem. We computed distances using the CS-2 code of Andrew Goldberg ([www.avglab.org/andrew](www.avglab.org/andrew)).

### *A.3. Using all affine projections*

The data analysis of Section 2 used projections into first and second order margins. The general theory developed later points to all affine projections as a natural base for analysis. In this section, we complete our analysis of the Plato data by looking at all affine projections.

In the following, $x$ and $z$ range over all binary 5-tuples. If $f(x)$ is the proportion of sentences in a fixed book (eg. *Republic*) with rhyming pattern $x$, the projection of $f$ in direction $z$ is

$$\sum_{x \cdot z = 0} f(x), \qquad \sum_{x \cdot z = 1} f(x).$$

TABLE 17
$d_V$ for Republic *to other books*

| Book | Vass. Dist. | Mass diff | Total | Rank |
|------|-------------|-----------|-------|------|
| Laws | 109 | 951 | 1060 | 5 |
| Phil. | 119 | 748 | 867 | 4 |
| Pol. | 112 | 952 | 1064 | 6 |
| Soph. | 82 | 97 | 179 | 1 |
| Tim. | 41 | 263 | 304 | 2 |
| Crit. | 71 | 675 | 746 | 3 |

Table 18

|       | (00010) | (01100) | (11000) |
|-------|---------|---------|---------|
| Rep.  | 1       | 2       | 1       |
| Tim.  | 2       | 3       | 2       |
| Soph. | 3       | 4       | 4       |
| Pol.  | 4       | 5       | 7       |
| Phil. | 7       | 6       | 3       |
| Laws. | 6       | 7       | 5       |
| Crit. | 5       | 1       | 6       |

To use the information that *Republic* was written early and *Laws* was written late, we find 5-tuples, $z$, that maximize

$$\left( \sum_{x \cdot z = 0} f(x) - \sum_{x \cdot z = 1} f(x) \right) - \left( \sum_{x \cdot z = 0} g(x) - \sum_{x \cdot z = 1} g(x) \right).$$

where $g(x)$ codes patterns for Laws. The largest three differences occur at $z = (00010), (01100)$ and $(11000)$. For each of these, we calculated

$$\sum_{x \cdot z = 0} h(x) - \sum_{x \cdot z = 1} h(x)$$

for each of the books (where $h$ codes the patterns for a particular book), and use the linear order of these values to order the books. The rank order resulting from the three binary 5-tuples, $z$, with the largest three differences above $z = (00010), (01100)$ and $(11000)$ are given in Table 18.

The first column thus gives the ranking: *Rep., Tim. Soph., Pol., Crit., Laws, Phil.* This is based on the difference between a single syllable (second from the end). It is close to, but not the same as the ranking based on adjusted second order margins found above. The other columns differ and show that not 'any old' projection gives the same ranking.

**Acknowledgments.** This paper is written in tribute to David Freedman with thanks for his integrity and brilliance.